\documentclass{amsproc}
\usepackage{graphicx}
\usepackage{amscd}
\usepackage{amsmath}
\usepackage{amsfonts}
\usepackage{amssymb}
\newtheorem{theorem}{Theorem}
\theoremstyle{plain}

\newtheorem{corollary}{Corollary}

\newtheorem{definition}{Definition}
\newtheorem{example}{Example}

\newtheorem{remark}{Remark}

\numberwithin{equation}{section}

\begin{document}
\title[The Amalgamation Property]{The Amalgamation Property in Classical Lebesgue-Riesz Spaces, Banach Spaces
with Almost Transitive Norm and Projection Constants}
\author{E.D. Positselskii}
\author{E.V. Tokarev}
\curraddr[E.V. Tokarev]{B.E. Ukrecolan, 33-81 Iskrinskaya str., 61005, Kharkiv-5, Ukraine}
\email[E.V. Tokarev]{tokarev@univer.kharkov.ua}
\thanks{This paper is in final form and no version of it will be submitted for
publication elsewhere.}
\dedicatory{Dedicated to the memory of S. Banach.}
\begin{abstract}In the paper it is shown that classes of finite equivalence that are generated
by Lebesgue-Riesz spaces $L_{p}$ have the amalgamation property if and only if
$p$ $\neq4,6,8,...$
\end{abstract}
\begin{abstract}This implies that the relative projection constant of a subspace
$A\hookrightarrow L_{p}\left[  0,1\right]  $ ($p$ $\neq4,6,8,...$) does not
depend on the position of $A$ in $L_{p}$. By the way the characterization of
2-dimensional subspaces of $l_{1}$ is obtained.
\end{abstract}
\begin{abstract}In the paper it is shown that classes of finite equivalence that are generated
by Lebesgue-Riesz spaces $L_{p}$ have the amalgamation property if and only if
$p$ $\neq4,6,8,...$
\end{abstract}
\begin{abstract}This implies that the relative projection constant of a subspace
$A\hookrightarrow L_{p}\left[  0,1\right]  $ ($p$ $\neq4,6,8,...$) does not
depend on the position of $A$ in $L_{p}$. By the way the characterization of
2-dimensional subspaces of $l_{1}$ is obtained.
\end{abstract}
\maketitle

\section{Introduction}

The amalgamation property for universal relational systems was introduced by
B. J\'{o}nsson in 1956 [1]. It was used to study partial ordered sets and
lattices, Boolean algebras and general algebraic relational systems.

V.I. Gurarii [2] independently discovered the same property in the isomorphic
theory of Banach spaces and showed that the set $\frak{M}$ of all
finite-dimensional Banach spaces enjoys it. This result (which is known under
the name \textit{the Gurarii's lemma on gluing of embeddings}) was used in [2]
to built the space $G$\ of almost universal disposition with respect to
$\frak{M}$.

In this paper we study the isometric version of the mentioned Gurarii's lemma
- i.e. the amalgamation property in the isometric theory of Banach spaces.

It will be shown that if a set $H(X)$ of all different finite dimensional
subspaces of a given Banach space $X$ (isometric subspaces in $H(X)$ are
identified) has the amalgamation property (all definitions will be given
later), then there exists a separable Banach space $G_{X}$, unique up to
almost isometry, which is finitely equivalent to $X$ and has properties that
are similar to those of the classical Gurarii space $G$:

\begin{itemize}
\item  Let $A$ and $B$ be isometric finite-dimensional subspaces of $X$;
$i:A\rightarrow B$ be the corresponding isometry. Then for every
$\varepsilon>0$ there exists an automorphism $u:X\rightarrow X$ such that
$\left\|  u\right\|  \left\|  u^{-1}\right\|  \leq1+\varepsilon$, which
extends $i$: $u\mid_{A}=i$. In particular, the norm of $G_{X}$ is almost transitive.

\item $G_{X}$ is an approximative envelope: every separable Banach space $Y$,
which is finitely representable in $G_{X}$ is $(1+\varepsilon)$-isomorphic to
a subspace of $G_{X}$ for every $\varepsilon>0$.
\end{itemize}

It will be convenient to regard the amalgamation property as a property of the
whole class $X^{f}$ of finite equivalence, generated by a given Banach space
$X$.

In this paper it will be shown that the class $\left(  L_{1}\right)  ^{f}$ has
the amalgamation property; that classes $\left(  L_{p}\right)  ^{f}$ does not
have the amalgamation property for $p=4$, $6$, $8$, ... (shortly: for
$p\in2\mathbb{N}\backslash\{2\}$). As it follows from [3], one W.~Rudin's
result [4] implies that classes $\left(  L_{p}\right)  ^{f}$ has the
amalgamation property for $p\in\left(  1,\infty\right)  \backslash\left(
2\mathbb{N}\backslash\{2\}\right)  $. The case $p=\infty$ is covered by the
mentioned above Gurarii's lemma. Worthy of note that in [4] was considered the
complex case. For real scalars the similar result was obtained in [5].

At the same time, in a general case the space $G_{X}$ (as above) is not of
almost universal disposition (with respect to the set $H(G_{X})$) in the sense
of [2]. This fact about $G_{L_{1}}$ may be deduced from one J. Pisier's result
[6]. The same result concerning $G_{L_{p}}$ will be obtained in this paper.

Some of results of this paper were announced in [7]; the tragical death of the
first named author terminated the publication of complete proofs.

\section{Definitions and notations}

\begin{definition}
Let $X$, $Y$ be Banach spaces. $X$ is \textit{finitely representable} in $Y$
(in symbols: $X<_{f}Y$) if for each $\varepsilon>0$ and for every finite
dimensional subspace $A$ of $X$ there exists a subspace $B$ of $Y$ and an
isomorphism $u:A\rightarrow B$ such that $\left\|  u\right\|  \left\|
u^{-1}\right\|  \leq1+\varepsilon$.

Spaces $X$ and $\ Y$ are said to be finitely equivalent, shortly: $X\sim_{f}%
Y$, if $X<_{f}Y$ and $Y<_{f}X$.

Any Banach space $X$ generates classes
\[
X^{f}=\{Y\in\mathcal{B}:X\sim_{f}Y\}\text{ \ and \ }X^{<f}=\{Y\in
\mathcal{B}:Y<_{f}X\}
\]
\end{definition}

For any two Banach spaces $X$, $Y$ their \textit{Banach-Mazur distance }is
given by
\[
d(X,Y)=\inf\{\left\|  u\right\|  \left\|  u^{-1}\right\|  :u:X\rightarrow Y\},
\]
where $u$ runs all isomorphisms between $X$ and $Y$ and is assumed, as usual,
that $\inf\varnothing=\infty$.

It is well known that $\log d(X,Y)$ defines a metric on each class of
isomorphic Banach spaces, where almost isometric Banach spaces are identified.

Recall that Banach spaces $X$ and $Y$ are \textit{almost isometric} if
$d(X,Y)=1$. Certainly, any almost isometric finite-dimensional Banach spaces
are isometric.

A set $\frak{M}_{n}$ of all $n$-dimensional Banach spaces, equipped with this
metric, is a compact metric space, called \textit{the Minkowski compact}
$\frak{M}_{n}$.

A disjoint union $\cup\{\frak{M}_{n}:n<\infty\}=\frak{M}$ is a separable
metric space, which is called the \textit{Minkowski space}.

Consider a Banach space $X$. Let $H\left(  X\right)  $ be a set of all its
different finite dimensional subspaces (isometric finite dimensional subspaces
of $X$ in $H\left(  X\right)  $ are identified). Thus, $H\left(  X\right)  $
may be regarded as a subset of $\frak{M}$, equipped with the restriction of
the metric topology of $\frak{M}$.

Of course, $H\left(  X\right)  $ need not to be a closed subset of $\frak{M}$.
Its closure in $\frak{M}$ will be denoted $\overline{H\left(  X\right)  }$.

From definitions\ it follows that $X<_{f}Y$ if and only if $\overline{H\left(
X\right)  }\subseteq\overline{H\left(  Y\right)  }$. Therefore, $X\sim_{f}Y$
if and only if $\overline{H\left(  X\right)  }=\overline{H\left(  Y\right)  }$.

There exists a one to one correspondence between classes of finite equivalence
$X^{f}$ and closed subsets of $\frak{M}$ of kind $\overline{H\left(  X\right)
}$.

Indeed, all spaces $Y$ from $X^{f}$ have the same set $\overline{H\left(
X\right)  }$. This set, uniquely determined by $X$ (or, equivalently, by
$X^{f}$), will be denoted by $\frak{M}(X^{f})$ and will be referred to as
\textit{the Minkowski's base of the class} $X^{f}$.

It will be convenient to introduce some new terminology.

A fifth $v=\left\langle A,B_{1},B_{2},i_{1},i_{2}\right\rangle $, where $A$,
$B_{1}$, $B_{2}\in\frak{M}(X^{f})$; $i_{1}:A\rightarrow B_{1}$ and$\ i_{2}%
:A\rightarrow B_{2}$ are isometric embeddings, will be called \textit{the }%
$V$\textit{-formation over} $\frak{M}(X^{f})$. The space $A$ will be called
\textit{the root of the }$V$\textit{-formation }$v$\textit{. }If there
exists\textit{\ }a triple $t=\left\langle j_{1},j_{2},F\right\rangle $, where
$F\in\frak{M}(X^{f})$; $\ j_{1}:B_{1}\rightarrow F$ and $j_{2}:B_{2}%
\rightarrow F$ are isometric embeddings such that $j_{1}\circ i_{1}=j_{2}\circ
i_{2}$, then the $V$-formation $v$ is said to be \textit{amalgamated in}
$\frak{M}(X^{f})$, and the triple $t$ is said to be its \textit{amalgam}.

Let $\operatorname{Amalg}(\frak{M}(X^{f}))$ be a set of all spaces
$A\in\frak{M}(X^{f})$ with the property:

\textit{Any }$V$\textit{-formation }$v$\textit{, which root is }%
$A$\textit{\ is amalgamated in }$\frak{M}(X^{f})$.

\begin{definition}
Let $X\in\mathcal{B}$ \ generates a class $X^{f}$\ with a Minkowski's base
$\frak{M}(X^{f})$. It will be said that $\frak{M}(X^{f})$ (and the class
$X^{f}$) has the amalgamation property if
\[
\frak{M}(X^{f})=\operatorname{Amalg}(\frak{M}(X^{f})).
\]
\end{definition}

\begin{definition}
For a Banach space $X$ its $l_{p}$-\textit{spectrum }$S(X)$ is given by
\[
S(X)=\{p\in\lbrack0,\infty]:l_{p}<_{f}X\}.
\]
\end{definition}

Certainly, if $X\sim_{f}Y$ then $S(X)=S(Y)$. So, the $l_{p}$-spectrum $S(X)$
may be regarded as a property of the whole class $X^{f}$.

\section{The amalgamation property in the class $\left(  L_{1}\right)  ^{f}$}

To show that a set $\frak{M}(\left(  l_{1}\right)  ^{f})$ has the amalgamation
property it will be needed some additional constructions.

\begin{definition}
Let $X$ be a Banach space with a basis $\left(  e_{n}\right)  $. This basis is
said to be Besselian if there exists such a constant $C>0$ that for any
$x=\sum x_{n}e_{n}\in X$%
\[
\left\|  x\right\|  \geq C\left(  \sum\left(  x_{n}\right)  ^{2}\right)
^{1/2}.
\]
\end{definition}

\begin{theorem}
Let $\left\langle X,\left(  e_{n}\right)  \right\rangle $ be a Banach space
with a Besselian bases $\left(  e_{n}\right)  $; $Y\hookrightarrow X$ be a
finite dimensional subspace of $X$ of dimension $\dim(Y)=m$.

There exists a sequence of points $Y_{0}=\left(  y_{n}\right)  _{n=1}^{\infty
}\subset\mathbb{R}^{m}$ such that $X$ may be represented as a Banach space
$\overline{X}$ of all rear-valued functions $f$, defined on $Y_{0}$ with a
norm
\[
\left\|  f\right\|  \overset{\operatorname{def}}{=}\left\|  f\right\|
_{\overline{X}}=\left\|  \sum f\left(  y_{i}\right)  e_{i}\right\|
_{X}<\infty;
\]
$Y$ may be represented as a subspace $\overline{Y}$ of $\overline{X}$ that
consists of all linear functions on $Y_{0}$ (equipped with the corresponding
restriction of the norm of $\overline{X}$). In other words, there exists an
isometry $T:X\rightarrow\overline{X}$ such that its restriction $T\mid
_{Y}:Y\rightarrow\overline{Y}$ is also an isometry.
\end{theorem}

\begin{proof}
Let $u=\sum u_{n}e_{n}$ and $v=\sum v_{n}e_{n}$ be elements of $X$. Since
$\left(  e_{n}\right)  $ is Besselian, it may be defined a scalar product
$\left\langle u,v\right\rangle =\sum u_{n}v_{n}$. Since $Y$ is of finite
dimension, there exists a projection $P:X\rightarrow Y$ such that
$\left\langle Px,y\right\rangle =\left\langle x,y\right\rangle $. Let
$y_{i}=Pe_{i}$.

Consider a linear space $\overline{X}$ of all formal sums $f=\sum f_{i}y_{i}$
and equip it with a norm
\[
\left\|  f\right\|  _{\overline{X}}=\left\|  \sum f_{i}y_{i}\right\|
_{\overline{X}}=\left\|  \sum f_{i}e_{i}\right\|  _{X}.
\]

Let a map $T:X\rightarrow\overline{X}$ be given by
\[
T(\sum u_{n}e_{n})=\sum u_{n}y_{n}\text{ \ for all }u=\sum u_{n}e_{n}\in X.
\]

For $u\in Y$ it may be computed its $i$'th coordinate
\[
\left(  Tu\right)  _{i}=u_{i}=\left\langle u,e_{i}\right\rangle =\left\langle
u,Pe_{i}\right\rangle =\left\langle u,y_{i}\right\rangle .
\]

Hence, $T$ maps $Y$ onto a linear space $\overline{Y}$ of all linear
functions, defined on $Y_{0}$.

So, $\overline{X}$, $\overline{Y}$, $Y_{0}$ and $T$ have desired properties.
\end{proof}

\begin{remark}
A set $Y_{0}=\left(  y_{i}\right)  \subset\mathbb{R}^{m}$ will be called an
incarnating set for a pair $Y\hookrightarrow X$. A pair $[\overline
{X},\overline{Y}]$ will be called an incarnation pair for $Y\hookrightarrow
X$. In what follows pairs $Y\hookrightarrow X$ and $[\overline{X},\overline
{Y}]$ may be identified.
\end{remark}

\begin{remark}
If a basis $\left(  e_{n}\right)  $ of $X$ is, in addition, a symmetric one, a
set $Y_{0}$ may be symmetrized: it may be considered instead $Y_{0}$ a central
symmetric set
\[
K=K\left(  Y\hookrightarrow X\right)  =Y_{0}\cup\left(  -Y_{0}\right)  =\{\pm
y\in\mathbb{R}^{m}:y\in Y_{0}\}
\]
as an incarnating set for a pair $Y\hookrightarrow X$ provided that
$\overline{X}\ $(resp., $\overline{Y}$) will be considered as the set of all
odd (resp., all odd linear) functions on $K$ with the corresponding norm.
\end{remark}

Consider a space $\mathbb{R}^{m}$. Let $\left(  e_{i}\right)  _{i=1}^{m}$be a
basis of $\mathbb{R}^{m}$ that will be assumed orthogonal (with respect to the
$l_{2}$-norm). Define on $\mathbb{R}^{m}$ an other norm, namely, put
\[
\left|  \left|  \left|  u\right|  \right|  \right|  =\left|  \left|  \left|
\sum\nolimits_{i=1}^{m}u_{n}e_{n}\right|  \right|  \right|  =\sum
\nolimits_{i=1}^{m}\left|  \left\langle u,e_{i}\right\rangle \right|
=\sum\nolimits_{i=1}^{m}\left|  u_{i}\right|  .
\]

\begin{theorem}
Let $m\in\mathbb{N}$, $K=\{\epsilon y_{i}:i\in\mathbb{N},$ $\epsilon
\in\{+,-\}\}\subset\mathbb{R}^{m}$ be a central symmetric set. $K$ is an
incarnating set for a pair $Y\hookrightarrow l_{1}$ where $Y$ is
$m$-dimensional subspace of $l_{1}$ if and only if $K$ is complete in
$\mathbb{R}^{m}$ and $\sum\{\left|  \left|  \left|  y_{i}\right|  \right|
\right|  :y_{i}\in K\}<\infty$.
\end{theorem}

\begin{proof}
Certainly, the completeness of $K$ in $\mathbb{R}^{m}$ is necessary.

The condition $\sum\{\left|  \left|  \left|  y_{i}\right|  \right|  \right|
:y_{i}\in K\}<\infty$ is also necessary. Indeed,
\begin{align*}
\sum\left|  \left|  \left|  y_{i}\right|  \right|  \right|   & =\sum
\nolimits_{y_{i}\in K}\sum\nolimits_{j=1}^{m}\left|  \left\langle y_{i}%
,e_{j}\right\rangle \right| \\
& =\sum\nolimits_{j=1}^{m}\sum\nolimits_{y_{i}\in K}\left|  \left\langle
y_{i},e_{j}\right\rangle \right|  =\sum\nolimits_{j=1}^{m}\left\|
e_{j}\right\|  =m<\infty.
\end{align*}

The sufficiency follows from a fact that a set $L(K)$ of all linear
functionals on $K$ may be embedded in $l_{1}$. Indeed, let $x\in L(K)$;
$x=\sum\{x_{i}y_{i}:y_{i}\in K\}$. Then
\begin{align*}
\left\|  x\right\|  _{l_{1}}  & =\sum\nolimits_{y_{i}\in K}\left|
\left\langle x,y_{i}\right\rangle \right|  =\sum\nolimits_{y_{i}\in K}\left|
\sum\nolimits_{y_{j}\in K}x_{j}\left\langle y_{i},y_{j}\right\rangle \right|
\\
& \leq\sum\nolimits_{y_{i}\in K}\max\left|  x_{j}\right|  \sum\nolimits_{y_{j}%
\in K}\left|  \left\langle y_{j},y_{i}\right\rangle \right|  =\max\left|
x_{j}\right|  \sum\nolimits_{y_{i}\in K}\left|  \left|  \left|  y_{i}\right|
\right|  \right|  .
\end{align*}
\end{proof}

\begin{remark}
Obviously, if linearly congruent sets $K$ and $K_{1}$ are incarnating sets for
pairs $Y\hookrightarrow l_{1}$ and $Z\hookrightarrow l_{1}$, then the linear
congruence $U:K\rightarrow K_{1}$ generates an isometric automorphism
$u:l_{1}\rightarrow l_{1}$, which restriction $u\mid_{Y}$ to $Y$ is an
isometry between $Y$ and $Z$.
\end{remark}

The following result shows how to reconstruct a pair $Y\hookrightarrow l_{1}$
by a given set $K=\{\pm y_{i}:i\in\mathbb{N}\}\subset\mathbb{R}^{m}$, which is
complete in $\mathbb{R}^{m}$ and satisfies the condition $\sum\{\left|
\left|  \left|  y_{i}\right|  \right|  \right|  :y_{i}\in K\}<\infty$.

Put
\[
\varsigma K=\cup\{\sum\{y_{i}\in K:\left\langle y_{i},u\right\rangle
>0\}:u\in\mathbb{R}^{m}\};\text{ \ }K^{\prime}=\overline{\operatorname{conv}%
\left(  \varsigma K\right)  },
\]
So, $K^{\prime}$ denotes the closure of a convex hull $\operatorname{conv}%
\left(  \varsigma K\right)  $ of $\varsigma K$.

\begin{theorem}
A central symmetric, complete in $\mathbb{R}^{m}$, set $K=\{\pm y_{i}%
:i\in\mathbb{N}\}\subset\mathbb{R}^{m}$ is an incarnating set for a pair
$Y\hookrightarrow l_{1}$ if and only if $K^{\prime}$ is congruent with the
unit ball
\[
B(Y^{\ast})\overset{\operatorname{def}}{=}\{y^{\prime}\in Y^{\ast}:\left\|
y^{\prime}\right\|  \leq1\}
\]
of the conjugate space $Y^{\ast}$.
\end{theorem}

\begin{proof}
\textbf{Necessity}.
\[
\left\|  x\right\|  _{Y}=\frac{1}{2}\sum\nolimits_{y_{i}\in K}\left|
\left\langle x,y_{i}\right\rangle \right|  =\max\{\left\langle
x,z\right\rangle :z\in K^{\prime}\}.
\]
Hence, $K^{\prime}=B(Y^{\ast})$.

\textbf{Sufficiency}. $K$ is an incarnating set for a some pair
$Z\hookrightarrow l_{1}$. Hence, as it was shown before, $K^{\prime}%
=B(Z^{\ast})$. However, by conditions of the theorem, $K^{\prime}=B(Y^{\ast})$.

Thus, $Y^{\ast}=Z^{\ast}$ and, hence, $Y=Z$.
\end{proof}

\begin{corollary}
There exists an one-to-one correspondence between pairs $Z\hookrightarrow
l_{1}$ and certain central symmetric complete subsets of $\mathbb{R}^{m}$.
\end{corollary}

The next result describe this correspondence in detail.

\begin{theorem}
Let $K$ be an incarnating set for a pair $Y\hookrightarrow l_{1}$. Vectors
$\left(  y_{i}\right)  $ of $K$ are parallel to edges of the unit ball
$B(Y^{\ast})$. A length of a sum of all vectors that are parallel to a given
edge $r$ of $B(Y^{\ast})$ is equal to the length of the edge $r$.
\end{theorem}

\begin{proof}
This result may be obtained as a consequence of [8]. Below it will be
presented its direct proof by induction on dimension.

If $\dim(K^{\prime})=1$ then $K^{\prime}=[-y,y]$ for some $y\in\mathbb{R} $;
$K=\{-y,y\}$ and the theorem is true trivially.

Assume that the theorem is true for $\dim(K^{\prime})=n$.

Let $K\in\mathbb{R}^{n+1}$; $\dim(K^{\prime})=n+1$.

An $n$-dimensional subspace $E\hookrightarrow\mathbb{R}^{n+1}$ is said to be a
support subspace (for $K^{\prime}$) if $E\cap K$ is complete in $E$ \ A set of
all support subspaces will be denoted by $E(K)$.

Let $E\in E(K)$; $r_{E}\in\mathbb{R}^{n+1}$; $r_{E}\neq0$. Let $r_{E}$ be
orthogonal to $E$. Let
\[
x_{E}=\sum\{y_{i}\in K:\left\langle y_{i},r_{E}\right\rangle >0\}.
\]

Then
\[
V_{E}=\left(  E\cap K\right)  ^{\prime}+x_{E}=\{z+x_{E}:z\in\left(  E\cap
K\right)  ^{\prime}\}
\]
is a facet of $K^{\prime}$.

Indeed, $\dim(K^{\prime})=n+1$; $\dim(V_{E})=n$; $V_{E}\subset K$. At the same
time $E+x_{E}$ is an affine manifold of the minimal dimension that contains
$V_{E}$. By the definition of $K^{\prime}$, $\left(  E+x_{E}\right)  \cap
K^{\prime}=V_{E}$. Hence, $V_{E}$ is one of maximal $n$-dimensional subsets of
$K^{\prime}$. Moreover,
\[
V_{E}=\{z\in K^{\prime}:\left\langle z,r_{E}\right\rangle =\left\langle
x_{E},r_{E}\right\rangle =\underset{x\in K^{\prime}}{\sup}\left\langle
x,r_{E}\right\rangle \}.
\]
So, to any support subspace $E\in E(K)$ corresponds a facet $V_{E}=\left(
E\cap K\right)  ^{\prime}+x_{E}$.

Conversely, let $V$ be a facet of $K^{\prime}$; $r_{V}$ be an external
perpendicular to $V$. Let
\begin{align*}
x_{V}  & =\sum\{y\in K:\left\langle y,r_{V}\right\rangle >0\};\\
E_{V}  & =\{x\in\mathbb{R}^{n+1}:\left\langle x,r_{V}\right\rangle =0\}.
\end{align*}
Then $V=\left(  E_{V}\cap K\right)  ^{\prime}+x_{V}$ and $E_{V}$ is a support subspace.

Let $E\in E(K)$. By the induction supposition, $E\cap K$ may be reconstructed
by the facet $V_{E}=\left(  E\cap K\right)  ^{\prime}+x_{V}$.

Since $K=\cup\{E\cap K:E\in E(K)\}$, then $K$ also may be reconstructed by
facets of $K^{\prime}$. Notice that every edge belongs to some facet of $K$.
\end{proof}

Now we are ready to prove the amalgamation property for $\frak{M}(\left(
l_{1}\right)  ^{f})$.

\begin{theorem}
The set $\frak{M}(\left(  l_{1}\right)  ^{f})$ has the amalgamation property.
\end{theorem}

\begin{proof}
We show that $H(l_{1})$ has the amalgamation property. Obviously, a closure
$\overline{H(l_{1})}=\frak{M}(\left(  l_{1}\right)  ^{f})$ also has this property.

Let $Y\hookrightarrow l_{1}$ and $Z\hookrightarrow l_{1}$. According to
preceding results these pairs determine corresponding incarnating sets
$K_{Y}=\left(  y_{i}\right)  $ and $K_{Z}=\left(  z_{j}\right)  $. Let $X$ be
isometrically embedded into $Y$ (by an operator $i:X\rightarrow Y)$ and into
$Z$ (by $j:X\rightarrow Z$).

Consider a conjugate $X^{\ast}$ and for every edge of its unit ball
$B(X^{\ast})$ choose subsets $\left(  y_{i}^{r}\right)  \subset K_{Y}=\left(
y_{i}\right)  $ and $\left(  z_{j}^{r}\right)  \subset K_{Z}=\left(
z_{j}\right)  $ that contains vectors parallel to $r$. Consider a set
$K_{XYZ}=\left(  x_{ij}\right)  =\{c_{ij}^{r}\}$, where $c_{ij}^{r}%
=y_{i}\left\|  z_{j}\right\|  /\left\|  r\right\|  $; $\left\|  r\right\|  $
denotes the length of the edge $r$. Clearly, $K_{XYZ}$\ is an incarnating set
for a subspace $W$ of $l_{1}$ that contains isometric images of both $Y$ and
$Z$ (say, $Y^{\prime}$ and $Z^{\prime}$) such that their intersection
$Y^{\prime}\cap Z^{\prime}=X^{\prime}$ is isometric to $X$; moreover,
corresponding isometries $u_{Y}:Y\rightarrow Y^{\prime}$ and $u_{Z}%
:Z\rightarrow Z^{\prime}$ transform embeddings $i:X\rightarrow Y$ and
$j:X\rightarrow Z$ to identical embeddings of $X^{\prime}$ into $Y^{\prime}$
and into $Z^{\prime}$ respectively. Clearly, this proves the theorem.
\end{proof}

\section{Almost $\omega$-homogeneous Banach spaces}

\begin{definition}
Let $X$ be a Banach space; $\mathcal{K}$ be a class of Banach spaces. $X$ is
said to be almost $\omega$-homogeneous with respect to $\mathcal{K}$ if for
any pair of spaces $A$, $B$ of $\mathcal{K}$ such that $A$ is a subspace of
$B$ ($A\hookrightarrow B$), every $\varepsilon>0$ and every isometric
embedding $i:A\rightarrow X$ there exists an isomorphic embedding $\hat
{\imath}:B\rightarrow X$, which extends $i$ (i.e., $\hat{\imath}|_{A}=i$) and
such, then
\[
\left\|  \hat{\imath}\right\|  \left\|  \hat{\imath}^{-1}\right\|
\leq(1+\varepsilon).
\]
\end{definition}

If $X$ is almost $\omega$-homogeneous with respect to $H(X)$ it will be simply
referred to as an almost $\omega$-homogeneous space.

\begin{theorem}
Any class $X^{f}$ of finite equivalence, whose Minkowski's set $\frak{M}%
(X^{f})$ has the amalgamation property contains a separable almost $\omega
$-homogeneous space $G_{X}$.

This space is unique up to almost isometry and is almost isotropic (in an
equivalent terminology, has an almost transitive norm).

This space is an approximative envelope of a class $X^{f}$: for any
$\varepsilon>0$ every separable Banach space which is finitely representable
in $X^{f}$ is $(1+\varepsilon)$-isomorphic to a subspace of $E_{X}$.
\end{theorem}

\begin{proof}
The proof of the theorem also literally repeats the Gurarii's one. We present
a sketch of proof for a sake of completeness (below it will be presented
another proof of this fact, based on a different idea).

Let us starting with any space $X_{0}$ of $X^{f}$. Consider a dense countable
subset $\{x_{i}:i<\infty;$ $\left\|  x_{i}\right\|  =1\}\subset X_{0}$. Any
finite subset $N$ of $\mathbb{N}$ defines a finite dimensional subspace
$U_{N}=\operatorname{span}\{x_{i}:i<\infty\}\hookrightarrow X_{0}$. Clearly,
there are only countable number of spaces of kind $U_{N}$. Identifying
isometric subspaces from
\[
\{U_{N}:N\subset\mathbb{N};\text{ \ }\operatorname{card}(N)<\infty\}
\]
it will be obtained a dense countable subset $\left(  Z_{i}\right)
_{i<\infty}\subset H(X_{0})$. Consider a set $\mathcal{F}$ of all triples
$\left\langle A,B,i\right\rangle $ where $A$, $B\in\frak{M}(X^{f})$ and there
exists an isometric embedding $i:A\rightarrow B$.

Let $n<m<\infty$ and\ let $\mathcal{F}\left(  n,m\right)  \subset\mathcal{F}$
be a subset of $\mathcal{F}$, which consists of such triples $\left\langle
A,B,i\right\rangle $ that $\dim\left(  A\right)  =n$; $\dim\left(  B\right)
=m $.

$\mathcal{F}\left(  n,m\right)  $ may be equipped with a metric
\[
\varrho\left(  \left\langle A,B,i\right\rangle ,\left\langle A_{1},B_{1}%
,i_{1}\right\rangle \right)  =\log d^{n,m}(\left\langle A,B,i\right\rangle
,\left\langle A_{1},B_{1},i_{1}\right\rangle ),
\]
where
\[
d^{n,m}(\left\langle A,B,i\right\rangle ,\left\langle A_{1},B_{1}%
,i_{1}\right\rangle )=\inf\{\left\|  u\right\|  \left\|  u^{-1}\right\|
:u:B\rightarrow B_{1};\text{ \ }i\circ u=i_{1}\}
\]
is a generalized Banach-Mazur distance. It is known [2] that $\left\langle
\mathcal{F}\left(  n,m\right)  ,\varrho\right\rangle $ is a compact metric
space. Hence, $\left\langle \mathcal{F},\varrho\right\rangle =\cup
_{n,m}\left\langle \mathcal{F}\left(  n,m\right)  ,\varrho\right\rangle $ is a
separable metric space and it may be chosen a countable dense subset
$\mathcal{F}^{0}$ of $\mathcal{F}$. Without loss of generality it may be
assumed that all spaces that are presented in triples from $\mathcal{F}^{0}$
are exactly those that belong to a defined before subset $\left(
Z_{i}\right)  \subset H(X_{0})$.

Using the amalgamation property for $\frak{M}(X^{f})$, for any pair
$A\hookrightarrow B\hookrightarrow X_{0}$ of subspaces of $X_{0}$ and any
isometric embedding $i:A\rightarrow C$, where $C\in\frak{M}(X^{f})$, it may be
constructed an extension of $X_{0}$, say, $X_{0}^{\prime}(A\hookrightarrow B;$
$i:A\rightarrow C)$ -- a separable Banach space that contains $X_{0}$ and
contains a triple $(A,B,C^{\prime})$, where $A\hookrightarrow B$;
$A\hookrightarrow C^{\prime}$ and a pair $A\hookrightarrow C^{\prime}$ is
isometric to the pair $iA\hookrightarrow C$ in a sense of the aforementioned
metric $\varrho$. It will be said that $X_{0}^{\prime}(A\hookrightarrow B;$
$i:A\rightarrow C)$ amalgamates over $X_{0}$ the $V$-formation
$(A\hookrightarrow B;$ $i:A\rightarrow C)$.

Now, proceed by induction.

Let $\left(  \frak{f}_{i}\right)  _{i<\infty}$ be a numeration of all triples
from $\mathcal{F}^{0}$. Construct a sequence of Banach spaces $\left(
X_{i}\right)  _{i<\infty}$ where it will be presented the space $X_{0}$ itself
as the first step of induction.

Let spaces $\left(  X_{i}\right)  _{i<\infty}$ \ be already constructed.

As $X_{n+1}$ it will be chosen the space $X_{n}^{\prime}(A\hookrightarrow B; $
$i:A\rightarrow C),$ where $A,B\in\left(  U_{N}\right)  _{N\subset\mathbb{N}}$;

$A\hookrightarrow B\hookrightarrow X_{n}$ be the $n$'th triple $\frak{f}_{n}$;

$\left(  i:A\rightarrow C\right)  =\frak{f}_{m+1}$, where $m$ is the least
number of a triple from $\mathcal{F}^{0}$ of kind $\left(  i:A\rightarrow
C\right)  $ (for a fixed $A$ and arbitrary $i$ and $C$) and be such that
$X_{n}$ does not amalgamate the $V$-formation $(A\hookrightarrow B;$
$i:A\rightarrow C)$.

Clearly, $X_{0}\hookrightarrow X_{1}\hookrightarrow X_{2}\hookrightarrow
..\hookrightarrow X_{n}\hookrightarrow...$. Let $X_{\infty}^{(1)}%
=\overline{\cup X_{i}}=\underset{\rightarrow}{\lim}X_{i}$. Now induction will
be continued, starting with $X_{\infty}^{(1)}$ instead of $X_{0}$. It will be
sequentially constructed spaces $X_{\infty}^{(2)}\hookrightarrow X_{\infty
}^{(3)}\hookrightarrow...\hookrightarrow X_{\infty}^{(n)}\hookrightarrow...$.

Their direct limit $\underset{\rightarrow}{\lim}X_{\infty}^{(n)}=G_{X}$ is the
desired space.

Proofs of the uniqueness of $G_{X}$ up to almost isometry and the property of
it to be almost isotropic are the same as in [2].
\end{proof}

\begin{definition}
($[9]$). Let $X$, $Y$ be Banach spaces; $Y\hookrightarrow X$. $Y$ is said to
be a reflecting subspace of $X$, shortly: $Y\prec_{u}X$, if for every
$\varepsilon>0$ and every finite dimensional subspace $A\hookrightarrow X$
there exists an isomorphic embedding $u:A\rightarrow Y$ such that $\left\|
u\right\|  \left\|  u^{-1}\right\|  \leq1+\varepsilon$, which is identical on
the intersection $Y\cap A$:
\[
u\mid_{Y\cap A}=Id_{Y\cap A}.
\]
\end{definition}

Clearly, $Y\prec_{u}X$ implies that $Y\sim_{f}X$.

\begin{definition}
($[10]$). A Banach space $E$ is said to be existentially closed in a class
$X^{f}$ if for any isometric embedding $i:E\rightarrow Z$ into an arbitrary
space $Z\in X^{f}$ its image $iE$ is a reflecting subspace of $Z$:
$iY\prec_{u}Z$.
\end{definition}

A class of all spaces $E$ that are existentially closed in $X^{f}$ is denoted
by $\mathcal{E}\left(  X^{f}\right)  $. In [10] it was shown that for any
Banach space $X$ the corresponding class $\mathcal{E}\left(  X^{f}\right)  $
is non-empty; moreover, any $Y<_{f}X^{f}$ may be isometrically embedded into
some $E\in\mathcal{E}\left(  X^{f}\right)  $ of the dimension $\dim
(E)=\max\{\dim(Y),\omega\}$.

\begin{theorem}
For any class $X^{f}$ such that $\frak{M}\left(  X^{f}\right)  $ has the
amalgamation property $E\in\mathcal{E}\left(  X^{f}\right)  $ if and only if
$E\in X^{f}$ and $E$ is almost $\omega$-homogeneous.
\end{theorem}

\begin{proof}
Let $E\hookrightarrow Z\in X^{f}$ and $E$ be almost $\omega$-homogeneous. Let
$A\hookrightarrow Z$ be a finite dimensional subspace; $E\cap A=B$ and
$\varepsilon>0$. Consider the identical embedding $id_{B}:B\rightarrow E$.
Since $B\hookrightarrow A$, $id_{B}$ may be extended to an embedding
$u:A\rightarrow E$ with $\left\|  u\right\|  \left\|  u^{-1}\right\|
\leq1+\varepsilon$. Thus, $E\in\mathcal{E}\left(  X^{f}\right)  $.

Conversely, let $E\in\mathcal{E}\left(  X^{f}\right)  $; $B\hookrightarrow E$;
$B\hookrightarrow A$ and $A\in\frak{M}\left(  X^{f}\right)  $.

Consider a space $Z$ such that $E\hookrightarrow Z$; $B\hookrightarrow
A\hookrightarrow Z$. Such space exists because of the amalgamation property of
$\frak{M}\left(  X^{f}\right)  $. Since $E\in\mathcal{E}\left(  X^{f}\right)
$, $E\prec_{u}Z$, i.e. there is an embedding $u:A\rightarrow E$ such that
$\left\|  u\right\|  \left\|  u^{-1}\right\|  \leq1+\varepsilon$, which is
identical on the intersection $E\cap A$: $u\mid_{E\cap A}=Id_{E\cap A}$. Since
$B\hookrightarrow A$ and $B\hookrightarrow E$, $B\hookrightarrow E\cap A$.
Clearly, $u$ extends the identical embedding $id_{B}$. Since $A\hookrightarrow
B$ and $\varepsilon$ are arbitrary, $E$ is almost $\omega$-homogeneous.
\end{proof}

\begin{remark}
This theorem gives an alternative proof of the first part of theorem 6
\end{remark}

\begin{corollary}
If $X^{f}$ is superreflexive and enjoys the amalgamation property then $G_{X}
$ is a norm one complemented subspace of any space $Z$ from the class $X^{f}$
that it contains.
\end{corollary}

\section{Projection constants}

Let $X$ be a Banach space; $Y\hookrightarrow X$. A \textit{projection
constant} $\lambda(A\hookrightarrow X)$ is given by
\[
\lambda(A\hookrightarrow X)=\inf\{\left\|  P\right\|  :\text{ }P:X\rightarrow
A;\text{ }P^{2}=P\},
\]
or, by words, $P$ runs all projections from $X$ onto $A$.

\textit{Relative projection constants} $\lambda(A,X)$ and $\lambda
(A,\mathcal{K})$, where $\mathcal{K}$ is a class of Banach spaces, are defined
as follows:
\[
\lambda(A,X)=\sup\{\lambda(iA\hookrightarrow X):\text{ \ }i:A\rightarrow X\},
\]
where $i$ runs all isometric embeddings of $A$ into $X$;
\[
\lambda(A,\mathcal{K})=\sup\{\lambda(A,X):X\in\mathcal{K}\}.
\]

The \textit{absolute projection constant }$\lambda(A)$ is just $\lambda
(A,\mathcal{B})$, where $\mathcal{B}$ denotes the class of all Banach spaces.

\begin{theorem}
Let $G$ be a separable almost $\omega$-homogeneous space; $A\hookrightarrow
G$. Then
\[
\lambda(A\hookrightarrow G)=\lambda(iA\hookrightarrow G)=\lambda(A,G),
\]
where $i$ is an arbitrary isometric embedding.
\end{theorem}

\begin{proof}
Obviously, any almost $\omega$-homogeneous space has the property:

\begin{itemize}
\item \textit{For any pair }$A$\textit{, }$B$\textit{\ of isometric subspaces
of }$G$\textit{\ under an isometry }$j:A\rightarrow B$\textit{\ and every
}$\varepsilon>0$\textit{\ there exists an isomorphic automorphism
}$u:G\rightarrow G$\textit{\ such that }$u\mid_{A}=j$\textit{\ and }$\left\|
u\right\|  \left\|  u^{-1}\right\|  \leq1+\varepsilon$\textit{.}
\end{itemize}

Clearly, this property implies the theorem.
\end{proof}

In [12] it was shown how to compute absolute projection constants of some
Banach spaces. Here a method, of computing relative projection constants will
be presented.

Let $X$ be an $N$-dimensional Banach space with a symmetric basis $\left(
e_{i}\right)  _{i=1}^{N}$; $Y$ be a subspace of $X$ of dimension $\dim(Y)=m$.
Let $K=K\left(  Y\hookrightarrow X\right)  =\left(  y_{i}\right)  _{i=1}%
^{2N}\subset\mathbb{R}^{m}$ be an incarnating (central symmetric) set for a
pair $Y\hookrightarrow X$. As $\overline{X}\ $(resp., as $\overline{Y}$) it
will be considered the set of all odd (resp., all odd linear) functions on $K$
with the corresponding norm. As in the previous section, a pair $\left[
\overline{X},\overline{Y}\right]  $ will be identified with $Y\hookrightarrow
X$.

Let $\frak{G}(X)$ be a group of all isometries of $X$; $\frak{G}_{0}$ be a
given group of (non degenerated) automorphisms of $\mathbb{R}^{m}$.

Assume that $K=\left(  y_{i}\right)  _{i=1}^{2N}\subset\mathbb{R}^{m}$ is
invariant under $\frak{G}_{0}$, i.e., $gK=K$ for all $g\in\frak{G}_{0}$. Then
$\frak{G}_{0}$ may be identified with a subgroup $\frak{G}_{0}^{\prime} $ of
$\frak{G}(X)$ under a group-isomorphism $\vartheta:\frak{G}_{0}\rightarrow
\frak{G}_{0}^{\prime}$: for every $g\in\frak{G}_{0}$ its image $g^{\prime
}=\vartheta g\in\frak{G}_{0}^{\prime}\subset\frak{G}(X)$ is given by
\[
\{g^{\prime}f\left(  y_{i}\right)  =f\left(  g^{-1}y_{i}\right)  :i=1\text{,
}2\text{, ..., }2N\}.
\]
In the future groups $\frak{G}_{0}$ and $\frak{G}_{0}^{\prime}$ will not be distinguished.

\begin{definition}
Let $\frak{G}$ be a subgroup of $\frak{G}(X)$. A projection $P:X\rightarrow Y
$ is said to be invariant under $\frak{G}$ if for all $g\in\frak{G}$ and $f\in
X$ it satisfies $P\left(  gf\right)  =gP\left(  f\right)  $.

A set of all projections that are invariant under $\frak{G}$ will be denoted
by $\operatorname{Inv}(\frak{G})$.
\end{definition}

Let $\frak{G}$ be a subgroup of the group $\mathbb{O}_{n}$\ of all orthogonal
matrices on $\mathbb{R}^{n}$. Let $m(g)$ be a normed Haar measure on
$\frak{G}$.

\begin{definition}
A subgroup $\frak{G}$ of the group $\mathbb{O}_{n}$ is said to be an ample
group if for all $u$, $v$, $x\in\mathbb{R}^{n}$
\[
\int\nolimits_{\frak{G}}gu(x,gv)dm(g)=n^{-1}(u,v)x,
\]
where $\left(  u,v\right)  $ is the usual scalar product.
\end{definition}

\begin{theorem}
Let $\frak{G}$ be a subgroup of the group $\mathbb{O}_{n}$; $\Gamma$ be a
subgroup of $\frak{G}$. If $\Gamma$ is an ample group then $\frak{G}$ is also ample.
\end{theorem}

\begin{proof}
Let $m(g)$ and $m^{\prime}(g)$ be normed Haar measures on $\frak{G}$ and on
$\Gamma$ respectively. Then
\begin{align*}
\int\nolimits_{\frak{G}}gu(x,gv)dm(g)  & =\int\nolimits_{\frak{G}}\left(
\int\nolimits_{\Gamma}g^{\prime}gu(x,g^{\prime}gv)dm^{\prime}(g)\right)
dm(g)\\
& =\int\nolimits_{\frak{G}}n^{-1}(gu,gv)xdm(g)=n^{-1}(u,v)x.
\end{align*}
\end{proof}

There may be presented a pair of important examples of ample groups.

\begin{example}
Let $\frak{G}_{1}$ be a subgroup of $\mathbb{O}_{n}$, which consists of all
operators that rearrange elements of a chosen basis $\left(  e_{i}\right)  $
and change their signs.

More exactly, let $\sigma$ be a rearrangement of $\left(  1,2,...,n\right)  $;
$\epsilon=\left(  \epsilon_{i}\right)  _{i=1}^{n}$; $\epsilon_{i}\in\{+,-\}$.
Then every element $g_{\epsilon}^{\sigma}$ of $\frak{G}_{1}$ acts on the
element $x=\sum\nolimits_{i=1}^{n}x_{i}e_{i}$ as follows:
\[
g_{\epsilon}^{\sigma}\left(  x\right)  =\sum\nolimits_{i=1}^{n}\epsilon
_{i}x_{i}e_{\sigma i}=u_{\epsilon}^{\sigma}.
\]
\end{example}

\begin{theorem}
$\frak{G}_{1}$ is an ample group.
\end{theorem}

\begin{proof}
Certainly,
\[
\sum\nolimits_{g\in\frak{G}_{1}}gu(x,gv)=\sum\nolimits_{\sigma,\epsilon
}u_{\epsilon}^{\sigma}\left(  x,v_{\epsilon}^{\sigma}\right)  =\sum
\nolimits_{\sigma}\sum\nolimits_{\epsilon}u_{\epsilon}^{\sigma}\left(
x,v_{\epsilon}^{\sigma}\right)  .
\]
Let
\[
\sum\nolimits_{\sigma,\epsilon}u_{\epsilon}^{\sigma}\left(  x,v_{\epsilon
}^{\sigma}\right)  =\sum\nolimits_{i=1}^{n}z_{i}e_{i}.
\]
Then
\begin{align*}
z_{i}  & =\sum\nolimits_{\epsilon}\epsilon_{i}u_{j_{i}}\left(  \epsilon
_{i}v_{j_{i}}x_{i}+\sum\nolimits_{l\neq i}\epsilon_{l}v_{j_{l}}x_{l}\right) \\
& =2^{n}u_{j_{i}}v_{j_{i}}x_{i}+u_{j_{i}}\sum\nolimits_{l\neq i}v_{j_{l}}%
x_{l}\left(  \sum\nolimits_{\epsilon}\epsilon_{j}\epsilon_{l}\right)  .
\end{align*}
Since $\sum\nolimits_{\epsilon}\epsilon_{j}\epsilon_{l}=0$ for $j\neq l$,
\[
z_{i}=2^{n}u_{j_{i}}v_{j_{i}}x_{i}.
\]
So,
\[
\sum\nolimits_{g\in\frak{G}_{1}}gu(x,gv)=2^{n}\sum\nolimits_{\sigma}%
\sum\nolimits_{i=1}^{n}u_{j_{i}}v_{j_{i}}x_{i}e_{i}=2^{n}\left(  n-1\right)
!\left(  u,v\right)  x.
\]
Since $\operatorname{card}\left(  \frak{G}_{1}\right)  =2^{n}n!$,
\[
\sum\nolimits_{g\in\frak{G}_{1}}gu(x,gv)=n^{-1}\operatorname{card}\left(
\frak{G}_{1}\right)  \left(  u,v\right)  x,
\]
i.e. $\frak{G}_{1}$ is an ample group.
\end{proof}

\begin{example}
Let the group $\frak{G}_{2}$ acts on a subspace of $\mathbb{R}^{n+1}$, which
is formed by such $x=\sum\nolimits_{i=1}^{m+1}x_{i}e_{i}\in\mathbb{R}^{n+1}%
$\ that $\sum\nolimits_{i=1}^{m+1}x_{i}=0$.

Every $g_{\pi}\in\frak{G}_{2}$ maps $x$ to $g_{\pi}x=y=\sum\nolimits_{i=1}%
^{m+1}x_{\pi i}e_{i}$ where $\pi$ is a rearrangement of $\left(
1,2,...,n+1\right)  $. So, $\operatorname{card}\left(  \frak{G}_{2}\right)
=\left(  n+1\right)  !$.
\end{example}

\begin{theorem}
$\frak{G}_{2}$ is an ample group.
\end{theorem}

\begin{proof}
Let $\sum\nolimits_{g\in\frak{G}_{2}}gu(x,gv)=\sum\nolimits_{i=1}^{n+1}%
z_{i}e_{i}$. Then
\begin{align*}
z_{j}  & =\sum\nolimits_{k=1}^{n+1}u_{k}\left(  n!x_{j}v_{k}+\sum
\nolimits_{i\neq j}x_{i}\left(  n-1\right)  !\left(  \sum\nolimits_{i\neq
k}v_{i}\right)  \right) \\
& =n!\left(  u,v\right)  x_{j}+\left(  n-1\right)  !\left(  \sum
\nolimits_{i\neq j}x_{i}\right)  \left(  \sum\nolimits_{k=1}^{n+1}\left(
u_{k}\sum\nolimits_{i\neq k}v_{i}\right)  \right)  .
\end{align*}
Since $\sum\nolimits_{i\neq j}x_{i}=-x_{j}$; $\sum\nolimits_{i\neq k}%
v_{i}=-v_{k}$,
\[
z_{j}=\left(  n!+\left(  n-1\right)  !\right)  \left(  u,v\right)
x_{j}=n^{-1}\operatorname{card}\left(  \frak{G}_{2}\right)  \left(
u,v\right)  x_{j}.
\]
Hence,
\[
\sum\nolimits_{g\in\frak{G}_{2}}gu(x,gv)=n^{-1}\operatorname{card}\left(
\frak{G}_{2}\right)  \left(  u,v\right)  x.
\]
\end{proof}

\begin{remark}
Certainly, the group $\frak{G}_{1}$ is a subgroup of the group of all
isometries of every $n$-dimensional Banach space with a symmetric basis
$\left(  e_{i}\right)  _{i=1}^{n}$.
\end{remark}

Let $E$ be a Banach space ; $\frak{G}(E)$ be the corresponding group of
isometries of $E$; $\left(  \frak{G}(E)\right)  ^{\#}$ be a set of all
(bounded, linear) automorphisms of $E$ that commute with all elements of
$\frak{G}(E)$.

As in [12] it will be said that $E$ is \textit{sufficiently symmetric} if
$\left(  \frak{G}(E)\right)  ^{\#}=\{\lambda Id_{E}\}_{\lambda\in\mathbb{R}}$.

\begin{theorem}
Let $E$ be a finite dimensional Banach space. $E$ is sufficiently symmetric if
and only if $\frak{G}(E)$ is an ample group.
\end{theorem}

\begin{proof}
Let $\frak{G}(E)$ be an ample group. Let $T:E\rightarrow E$ commutes with all
elements of $\frak{G}(E)$. Then
\begin{align*}
n^{-1}\left(  u,Tv\right)  x  & =\int\nolimits_{\frak{G}(E)}gu\left(
gTv,x\right)  dm(g)\\
& =\int\nolimits_{\frak{G}(E)}gu\left(  T^{\ast}x,gv\right)  dm(g)=n^{-1}%
\left(  u,v\right)  T^{\ast}x.
\end{align*}
Put $u=v$. Then $T^{\ast}x=\left(  v,Tv\right)  \left\|  v\right\|  ^{-2}x$
for all $x\in E$. Hence, $T^{\ast}x=\lambda x$ and $T^{\ast}=\lambda
Id_{E^{\ast}}$.

Conversely, let $E$ be sufficiently symmetric. Fix $u$, $v\in E$ and put for
$x\in E$%
\[
Tx=\int\nolimits_{\frak{G}(E)}gu\left(  x,gv\right)  dm(g).
\]

Let $g_{1}\in\frak{G}(E)$. Then
\begin{align*}
Tg_{1}x  & =\int\nolimits_{\frak{G}(E)}gu\left(  g_{1}x,gv\right)
dm(g)=g_{1}\int\nolimits_{\frak{G}(E)}g_{1}^{-1}gu\left(  x,g_{1}%
^{-1}gv\right)  dm(g)\\
& =g_{1}\int\nolimits_{\frak{G}(E)}g_{1}^{-1}gu\left(  x,g_{1}^{-1}gv\right)
dm(g_{1}^{-1}g)=g_{1}Tx.
\end{align*}

Hence, there exists a constant $c=c\left(  x,y\right)  $ such that $T=c\left(
x,y\right)  Id_{E}$. Indeed,
\[
\int\nolimits_{\frak{G}(E)}\left(  gu,y\right)  \left(  x,gv\right)
dm(g)=\int\nolimits_{\frak{G}(E)}\left(  u,gy\right)  \left(  gx,v\right)
dm(g)=c\left(  x,y\right)  \left(  u,v\right)  .
\]
At the same time,
\[
\int\nolimits_{\frak{G}(E)}\left(  gu,y\right)  \left(  x,gv\right)
dm(g)=\left(  Tx,y\right)  =c\left(  u,v\right)  \left(  x,y\right)  .
\]
Obviously, $c(x,y)=\left(  \dim(E)\right)  ^{-1}=n^{-1}$. So, $\int
\nolimits_{\frak{G}(E)}gu\left(  x,gv\right)  dm(g)=n^{-1}\left(  u,v\right)
x$.

Hence, $\frak{G}(E)$ is an ample group.
\end{proof}

Now it will be presented a formula for computing some relative projection
constants. Let $Y\hookrightarrow X$ be realized as an incarnation pair
$[\overline{X},\overline{Y}]$ with an incarnating central symmetric set
$K=\left(  y_{i}\right)  _{i=1}^{2N}\subset\mathbb{R}^{n}$, invariant under an
ample group $\frak{G}$.

Let $\mathcal{M}$ be a set of all mappings $\mu:K\rightarrow\mathbb{R}^{n}$
such that $\mu\left(  gu\right)  =g\mu\left(  u\right)  $ for any
$g\in\frak{G}$, $u\in K$.

\begin{theorem}
Let $P:X\rightarrow Y$ be a linear mapping. $P\in\operatorname{Inv}(\frak{G})$
is and only if $P$ has a representation
\[
P(f)=n\left(  \sum\nolimits_{u\in K}\left(  u,\mu\left(  u\right)  \right)
\right)  ^{-1}\sum\nolimits_{u\in K}f\left(  u\right)  \mu\left(  u\right)  .
\]
\end{theorem}

\begin{proof}
Certainly, if $P$ has such representation then $P$ is a linear map.

Put $c=n\left(  \sum\nolimits_{u\in K}\left(  u,\mu\left(  u\right)  \right)
\right)  ^{-1}$. Let $m(g)$ be a normed Haar measure on $\frak{G}$. Then for
$y\in Y$
\begin{align*}
Py  & =\sum\nolimits_{u\in K}\left(  u,y\right)  \mu\left(  u\right)
=c\sum\nolimits_{u\in K}\int\nolimits_{\frak{G}}\left(  gu,y\right)
\mu\left(  gu\right)  dm\left(  g\right) \\
& =c\sum\nolimits_{u\in K}\int\nolimits_{\frak{G}}\left(  gu,y\right)
\mu\left(  u\right)  dm\left(  g\right)  =n^{-1}c\sum\nolimits_{u\in K}\left(
u,\mu\left(  u\right)  \right)  y=y.
\end{align*}

Let us show that $P$ is invariant under $\frak{G}$.
\begin{align*}
P\left(  gf\right)   & =P\left(  fg^{-1}\right)  =c\sum\nolimits_{u\in
K}f\left(  g^{-1}\mu\right)  \mu\left(  u\right)  =c\sum\nolimits_{v\in
K}f\left(  v\right)  \mu\left(  gv\right) \\
& =g\left(  c\sum\nolimits_{v\in K}f\left(  v\right)  \mu\left(  v\right)
\right)  =gP(f).
\end{align*}

To show the necessity of the conditions of the theorem assume that
$P\in\operatorname{Inv}(\frak{G})$. Let $\{e_{1},e_{2},...,e_{n}\}$ be a basis
of $\mathbb{R}^{n}$.

Consider linear functionals $\left(  \varphi_{i}\right)  _{i=1}^{n}$ over $X$,
which are given by
\[
\varphi_{i}\left(  f\right)  =\left(  Pf,e_{i}\right)  =\sum\nolimits_{v\in
K}f\left(  u\right)  \mu_{i}\left(  u\right)  ,
\]
where $\mu_{i}\in X$. In a vector form this equality may be written as
\[
P\left(  f\right)  =\sum\nolimits_{u\in K}f\left(  u\right)  \mu\left(
u\right)  ;\text{ \ }\mu=\{\mu_{1},\mu_{2},...,\mu_{n}\}.
\]

Let $g\in\frak{G}$. Then
\[
gP\left(  f\right)  =P\left(  gf\right)  =P\left(  fg^{-1}\right)
=\sum\nolimits_{u\in K}f\left(  g^{-1}u\right)  \mu\left(  u\right)
=\sum\nolimits_{u\in K}f\left(  u\right)  \mu\left(  gu\right)  .
\]
Hence,
\[
\sum\nolimits_{u\in K}f\left(  u\right)  g\mu\left(  u\right)  =\sum
\nolimits_{u\in K}f\left(  u\right)  \mu\left(  gu\right)
\]
for all $f\in X$. Since $X$ consists of odd functions, $\mu\left(  -u\right)
=\mu\left(  u\right)  $. Hence, $\mu\in\mathcal{M}$.

Let $y\in Y$. Since $\mu\left(  gu\right)  =g\mu\left(  u\right)  $ and
$\frak{G}$ is an ample group,
\begin{align*}
P(y)  & =\sum\nolimits_{u\in K}\left(  u,y\right)  \mu\left(  u\right)
=\sum\nolimits_{u\in K}\int\nolimits_{\frak{G}}\left(  y,gu\right)  \mu\left(
gu\right)  dm\left(  g\right)  =\\
& =\sum\nolimits_{u\in K}\int\nolimits_{\frak{G}}\left(  y,gu\right)
g\mu\left(  u\right)  dm\left(  g\right)  =n^{-1}\sum\nolimits_{u\in K}\left(
u,\mu(u)\right)  y.
\end{align*}

Hence, $n^{-1}\sum\nolimits_{u\in K}\left(  u,\mu(u)\right)  =1$ and,
consequently,
\[
P(f)=n\left(  \sum\nolimits_{u\in K}\left(  u,\mu\left(  u\right)  \right)
\right)  ^{-1}\sum\nolimits_{u\in K}f\left(  u\right)  \mu\left(  u\right)  .
\]
\end{proof}

An easy consequence of this result is a formula for computing projection constants.

\begin{corollary}
Let $Y\hookrightarrow X$ be realized as an incarnation pair $[\overline
{X},\overline{Y}]$ with an incarnating central symmetric set $K=\left(
y_{i}\right)  _{i=1}^{2N}\subset\mathbb{R}^{n}$, which is invariant under an
ample group $\frak{G}$. Let $\mathcal{M}$ be a set of all mappings
$\mu:K\rightarrow\mathbb{R}^{n}$ such that $\mu\left(  gu\right)  =g\mu\left(
u\right)  $ for any $g\in\frak{G}$, $u\in K$. Let $P:X\rightarrow Y$ be a
linear projection.

Then the projection constant $\lambda\left(  Y\hookrightarrow X\right)  $ is
equal to
\[
\lambda\left(  Y\hookrightarrow X\right)  =\inf_{\mu\in\mathcal{M}}\{n\left(
\sum\nolimits_{u\in K}\left(  u,\mu\left(  u\right)  \right)  \right)
^{-1}\}\sup_{f\in X}\{\sum\nolimits_{u\in K}f\left(  u\right)  \mu\left(
u\right)  \}.
\]
\end{corollary}

\begin{proof}
It is known (see e.g. [13]) that if $K$ is invariant under $\frak{G}$ then
\[
\lambda\left(  Y\hookrightarrow X\right)  =\inf\{\left\|  P\right\|
:P\in\operatorname{Inv}(\frak{G});\text{ \ }P:X\rightarrow Y;\text{ \ }%
P^{2}=P\}.
\]
So, the desired result follows from the representation of $P$ as above.
\end{proof}

\section{Amalgamation property in classes $\left(  L_{2n}\right)  ^{f}$ for a
natural $n>1$}

\begin{theorem}
If $p\in2\mathbb{N}\backslash\{2\}=\{4,6,8,...\}$ then the set $\frak{M}%
(\left(  L_{p}\right)  ^{f})$ does not have the amalgamation property.
\end{theorem}

\begin{proof}
Assume that $\frak{M}(\left(  L_{p}\right)  ^{f})$ has the amalgamation
property. Then, by the theorem 6, there exists a separable almost $\omega
$-homogeneous space $W$, which belongs to $\left(  L_{p}\right)  ^{f}$. By the
theorem 7, $W$ must be existentially closed in $\left(  L_{p}\right)  ^{f}$.
According to [11], $L_{p}$ is an envelope of $\left(  L_{p}\right)  ^{f}$ and,
hence, contains a subspace that is isomorphic to $W$. By the corollary 2, $W$
is 1-complemented subspace of $L_{p}$. From the theorem 2 it follows that
$L_{p}$ is almost isometric to a complemented subspace of $W$. Since every
isometric image of $L_{p}$ in $L_{p}$ is orthogonal complemented, it is clear
that for any pair $A$, $B$ of isometric finite dimensional subspaces of
$L_{p}$ their projection constants must be equal: $\lambda\left(
A\hookrightarrow L_{p}\right)  =\lambda\left(  B\hookrightarrow L_{p}\right)
$. Let us show that if $p\in2\mathbb{N}\backslash\{2\}$ then $L_{p}$ contains
a pair of Euclidean 2-dimensional subspaces that have different projection
constants. Notice, that it is enough to find such a pair in some
$l_{p}^{\left(  n\right)  }$.

Let $p=2n$; $n\geq2$. Let $K=(y_{k})_{k=1}^{2n+2}$ be a set of vertices of
symmetric polygon. Let $(y_{k})_{k=1}^{2n+2}$ has polar coordinates
\[
\{\left(  \varphi_{k},R_{k}\right)  =\left(  \pi k/\left(  n+1\right)
;1\right)  :k=0,1,...,2n+1\}.
\]

Let $L(K)$ be a set of all odd linear functions on $K$. The norm of an element
$u=\left(  u_{1},u_{2}\right)  \in L(K)$ in the space $l_{2n}^{\left(
2n+2\right)  }$ (that is considered as a space of all odd functions on $K$
with the corresponding norm) is calculated by
\begin{align*}
\left\|  u\right\|  ^{2n}  & =\sum\nolimits_{n=0}^{2n+1}\left(  u_{1}%
\cos\left(  \frac{k\pi}{n+1}\right)  +u_{2}\sin\left(  \frac{k\pi}%
{n+1}\right)  \right)  ^{2n}\\
& =\left(  \left\|  u\right\|  _{2}\right)  ^{2n}\sum\nolimits_{n=0}%
^{2n+1}\left(  \cos\left(  \frac{k\pi}{n+1}\right)  \sin\phi+\sin\left(
\frac{k\pi}{n+1}\right)  \cos\phi\right)  ^{2n}\\
& =\left(  \left\|  u\right\|  _{2}\right)  ^{2n}\sum\nolimits_{n=0}%
^{2n+1}\left[  \sin\left(  \frac{k\pi}{n+1}+\phi\right)  \right]  ^{2n},
\end{align*}
where $\left\|  u\right\|  _{2}$ $=\sqrt{\left|  u_{1}\right|  ^{2}+\left|
u_{2}\right|  ^{2}}$. Regard $\left\|  u\right\|  ^{2n}$ as an function of
$\phi$; $\left\|  u\right\|  ^{2n}=F(\phi)$. Immediately, the derivative
$F^{\prime}(\phi)=0$. So, the sum is a constant $c$ and $\left\|  u\right\|
=c\left\|  u\right\|  _{2}$, i.e. $L(K)$ is isometric to the 2-dimensional
Euclidean space.

Let $\Delta_{\varphi}$ be a linear operator of rotation of the plane in the
positive direction by the angle $\varphi$. Put $K(\varphi)=\Delta_{\varphi
}\left(  K\right)  $; $K^{\left[  m\right]  }=\cup_{i=0}^{m-1}K(\pi i/nm)$.
Then
\[
\left\|  u\right\|  _{L(K^{\left[  m\right]  })}^{2n}=\sum\nolimits_{i=0}%
^{m-1}\left\|  \Delta_{\pi i/nm}\left(  u\right)  \right\|  _{L(K^{\left[
m\right]  })}^{2n}=m\left\|  u\right\|  _{L(K)}^{2n}.
\]

Hence 2-dimensional subspace $L(K^{\left[  m\right]  })$ of $l_{2n}^{\left(
\left(  2n+2\right)  m\right)  }$ is also isometric to the Euclidean one.

Let us show that
\[
\lambda(L\left(  K^{\left[  m\right]  }\right)  \hookrightarrow l_{2n}%
^{((2n+2)m)}(K^{\left[  m\right]  }))<\lambda(L\left(  K\right)
\hookrightarrow l_{2n}^{(2n+2)}\left(  K\right)  ).
\]

Certainly, this would imply the needed result.

For simplicity of calculations assume that $n=2$; the result in a general case
is obtained in a similar way.

Notice that $K$ is invariant under the apple group $\frak{G}_{1}$. Hence, it
may be used the corollary 3.
\begin{align*}
\lambda(L(K^{\left[  m\right]  })  & \hookrightarrow l_{4}^{(6m)}(K^{\left[
m\right]  }))=\\
& =2\left[  \sum\nolimits_{u\in K^{\left[  m\right]  }}\left(  u,u\right)
\right]  ^{-1}\sup_{f\in l_{4}^{(6m)}\left(  K^{\left[  m\right]  }\right)
}\{\left\|  f\right\|  ^{-1}\left\|  \sum\nolimits_{u\in K^{\left[  m\right]
}}\left(  u\right)  fu\right\|  _{l_{4}}\}\\
& =\frac{6}{n}\left\|  A_{n}\right\|  _{l_{4}\left(  K\right)  \rightarrow
L(K)},\text{ \ \ \ \ where }A_{n}\left(  f\right)  =\sum\nolimits_{u\in
K^{\left[  m\right]  }}f(u)u.
\end{align*}

The adjoint operator $A_{n}^{\ast}:L^{\ast}(K^{\left[  m\right]  })\rightarrow
l_{4/3}\left(  K^{\left[  m\right]  }\right)  $ is given by $A_{n}^{\ast
}\left(  v\right)  =g$, where $g\left(  y\right)  =\left(  v,y\right)  $;
$y\in K^{\left[  m\right]  }$. Since $\left\|  A_{n}\right\|  =\left\|
A_{n}^{\ast}\right\|  $,
\begin{align*}
\lambda(L(K^{\left[  m\right]  })  & \hookrightarrow l_{4}^{(6m)}(K^{\left[
m\right]  }))\\
& =\frac{1}{3n}\sup\{\left\|  v\right\|  ^{-1}\left\|  v\right\|
_{l_{3/4}\left(  K^{\left[  m\right]  }\right)  }:v\in L^{\ast}(K^{\left[
m\right]  })\}\\
& =\frac{1}{3n}\left(  \frac{3\sqrt{n}}{2\sqrt{2}}\right)  ^{1/2}\sup
_{v}\{\left(  v,v\right)  ^{-1/2}\left(  \sum\nolimits_{i=0}^{m-1}\left\|
\Delta_{\pi i/3n}\left(  v\right)  \right\|  _{l_{4/3}(K^{\left[  m\right]
})}^{4/3}\right)  ^{3/4}\\
& \leq\frac{1}{3n}\left(  \frac{3\sqrt{n}}{2\sqrt{2}}\right)  ^{1/2}\left(
\sum\nolimits_{i=0}^{m-1}\sup_{v}\left[  \left(  v,v\right)  ^{-1/2}\left\|
\Delta_{\pi i/3n}\left(  v\right)  \right\|  _{l_{4/3}(K^{\left[  m\right]
})}\right]  ^{4/3}\right)  ^{3/4}\\
& =\left(  \frac{1}{6\sqrt{2}}\right)  ^{1/2}\sup_{v}\left[  \left(
v,v\right)  ^{-1/2}\left\|  v\right\|  _{l_{3/4}\left(  K^{\left[  m\right]
}\right)  }\right]  =\lambda(L\left(  K\right)  \hookrightarrow l_{2n}%
^{(2n+2)}\left(  K\right)  ).
\end{align*}

The equality may take place only if the exact upper bound in the expression
$\sup_{v}\left[  \left(  v,v\right)  ^{-1/2}\left\|  \Delta_{\pi i/3n}\left(
v\right)  \right\|  _{l_{4/3}(K^{\left[  m\right]  })}\right]  $ is attained
for all $i=0$, $1$, ..., $n-1$ on the same vector $v$ .

However it is easy to check that
\[
\left\|  \left(  1/2,\sqrt{3}/2\right)  \right\|  _{l_{4/3}\left(  K\right)
}>\left\|  \left(  \sqrt{2}/2,\sqrt{2}/2\right)  \right\|  _{l_{4/3}\left(
K\right)  }.
\]

Hence, if $n$ is sufficiently large, it will take place the strong inequality
\[
\lambda(L(K^{\left[  m\right]  })\hookrightarrow l_{4}^{(6m)}(K^{\left[
m\right]  }))<\lambda(L\left(  K\right)  \hookrightarrow l_{4}^{(6)}\left(
K\right)  )
\]
and, consequently, the set $\frak{M}(\left(  L_{p}\right)  ^{f})$ does not
have the amalgamation property.
\end{proof}

\section{Properties of spaces $L_{1}\left[  0,1\right]  $ and $l_{1}$}

\begin{theorem}
Space $L_{1}\left[  0,1\right]  $ is almost $\omega$-homogeneous.
\end{theorem}

\begin{proof}
Since $\left(  G_{X}\right)  ^{\ast\ast}$ is a complemented subspace of
$\left(  L_{1}\left[  0,1\right]  \right)  ^{\ast\ast}$, $G_{X}$ is an
$\mathcal{L}_{1,1+0}$-space in a sense of J. Lindenstrauss and A. Pelczynski
[11] and, hence, is a separable $L_{1}\left(  \mu\right)  $-space. All
pairwice non-isometric separable $L_{1}\left(  \mu\right)  $-spaces may be
listed. They are:
\begin{align*}
& l_{1};\text{ }L_{1}\left[  0,1\right]  \text{\ \ }l_{1}\oplus_{1}%
L_{1}\left[  0,1\right]  ,\\
\{l_{1}^{(n)}\oplus_{1}L_{1}\left[  0,1\right]   & :n\in\mathbb{N}\}\text{ and
}l_{1}\left(  L_{1}\left[  0,1\right]  \right)  \overset{\operatorname{def}%
}{=}\left(  \sum\oplus L_{1}\left[  0,1\right]  \right)  _{l_{l}}%
\end{align*}

It is easy to check that only the space $L_{1}\left[  0,1\right]  $ is
existentially closed. Any other space cannot be isometrically embedded in
$L_{1}\left[  0,1\right]  $ in a such way that its image will be a reflecting
subspace of $L_{1}\left[  0,1\right]  $.
\end{proof}

\begin{corollary}
For any isometric finite dimensional subspaces $A$, $B$ of $L_{1}\left[
0,1\right]  $ their relative projection constants are equal:
\[
\lambda\left(  A\hookrightarrow L_{1}\left[  0,1\right]  \right)
=\lambda\left(  B\hookrightarrow L_{1}\left[  0,1\right]  \right)  .
\]
\end{corollary}

\begin{proof}
This follows from the preceding theorem and the theorem 8.
\end{proof}

The next result describes two-dimensional subspaces of $L_{1}\left[
0,1\right]  $. Before its formulation recall that a point $z$ of the unit
sphere $\partial B(X)=\{z\in X:\left\|  z\right\|  =1\}$ of a Banach space $X$
is said to be \textit{extreme} if for any pair $z_{1}$, $z_{2}\in\partial
B(X)$ the condition $\left(  z_{1}+z_{2}\right)  /2=z$ implies that
$z_{1}=z_{2}=z$.

The set of all extreme points of $\partial B(X)$ is denoted by
$\operatorname{ext}(X) $.

\begin{theorem}
Let $Z$ be 2-dimensional Banach space. The following conditions are equivalent:

\begin{enumerate}
\item $Z^{\ast}$ is isometric to a 2-dimensional subspace of $l_{1}$.

\item $\operatorname{ext}(Z)$ is of linear zero measure $\operatorname{mes}%
\operatorname{ext}(Z)$ on the boundary $\partial B(Z)$ of the unit cell of $Z$.
\end{enumerate}
\end{theorem}

\begin{proof}
$\left(  2\Rightarrow1\right)  $. Let $V=\operatorname{ext}\left(  Z\right)  $.

Since $\operatorname{mes}(V)=0$ on $\partial B(Z)$, and $V$ is a closed subset
of $\partial B(Z)$, $V$ is nowhere dense in $\partial B(Z)$.

Hence, $\partial B(Z)$ is a closure of the union of edges of $\partial B(Z)$.
Let $\{r_{i}:i\in I\}$ be their numeration. Let $\left|  r_{i}\right|  $ be
the length of the edge $r_{i}$.

Since $\sum\nolimits_{i\in I}\left|  r_{i}\right|  =\operatorname{mes}\left(
\partial B\left(  Z\right)  \right)  <\infty$, the set $I$ is at most countable.

Consider a edge $r_{i}$ and denote its ends by $x_{2i-1}$ and $x_{2i}$ (the
numeration of this ends is formed by the positive direction of path-tracing on
$\partial B(Z)$).

Let $K=\{y_{i}=\left(  x_{2i}-x_{2i-1}\right)  /2:i<\infty\}$. Let us show
that $K$ is an incarnation set for a pair $Z^{\ast}\hookrightarrow l_{1}$. It
is sufficient to show that $K^{\prime}=B(Z)$ (notations are as in the section 3).

Let $x\in V$. Let $R_{x}$ be a set of all edges of the unit circle $B(Z)$ that
are appeared on the way (along the unit circumference $\partial B\left(
Z\right)  $) in the positive direction from $-x$ up to $x$. It will be shown
that $\sum\{y_{i}:r_{i}\in R_{x}\}=x$.

Fix $\varepsilon>0$. Since $V$ is a metric compact of measure $0$, there
exists a finite number of pairwice disjoint arcs of $\partial B(Z)$ that cover
$V$ and whose total length is less then $\varepsilon$. Let $z_{1}$, $z_{2}$,
..,$z_{N}$ be ends of these arcs ( in the positive direction of the path-tracing).

Certainly, $\left(  z_{1}+x\right)  /2+\left(  z_{2}-z_{1}\right)
/2+...+\left(  x-z_{N}\right)  /2=x$ and
\[
\left\|  \left(  z_{1}+x\right)  /2\right\|  +\left\|  \left(  z_{2}%
-z_{1}\right)  /2\right\|  +...+\left\|  \left(  x-z_{N}\right)  /2\right\|
\leq\varepsilon.
\]

By the construction,
\[
\left\|  \sum\{y_{i}:r_{i}\in R_{x}\}-\left[  \left(  z_{1}+x\right)
/2+\left(  z_{2}-z_{1}\right)  /2+...+\left(  x-z_{N}\right)  /2\right]
\right\|  \leq\varepsilon.
\]
Hence,
\[
\left\|  x-\sum\{y_{i}:r_{i}\in R_{x}\}\right\|  \leq2\varepsilon.
\]
Since $\varepsilon$ is arbitrary, $x=\sum\{y_{i}:r_{i}\in R_{x}\}$.

Recall that $x\in\partial B\left(  Z\right)  $ is said to be an
\textit{exposed point} if there exists such $f\in Z^{\ast}$ that $f\left(
x^{\prime}\right)  <f\left(  x\right)  $ for all $x^{\prime}\neq x$,
$x^{\prime}\in\partial B\left(  Z\right)  $.

Let $\Pi$ be a supporting line for $B\left(  Z\right)  $ such that $\Pi\cap
B\left(  Z\right)  =\{x\}$.

Let $u$ be such that $\left(  u,y\right)  >0$ for all $y\in\Pi$. Then
\[
\{y_{i}\in K:\left(  y_{i},u\right)  >0\}=\{y_{i}:r_{i}\in R_{x}\}.
\]
Hence,
\[
x=\sum\{y_{i}:r_{i}\in R_{x}\}=\sum\{y_{i}:y_{i}\in K;\text{ }\left(
y_{i},u\right)  >0\}\in K^{\prime}.
\]

According to the known S. Straszewicz theorem [14], the set of all exposed
points of $B(Z)$ is dense in $\operatorname{ext}(Z)=V$. So, $V\subset
K^{\prime}$ and, consequently, $B(Z)\subset K^{\prime}$.

Conversely, let $x\in K^{\prime}$; $x=\sum\{y_{i}:y_{i}\in K;$ $\left(
y_{i},u\right)  >0\}$. Let $\Pi$ be a supporting line for $B\left(  Z\right)
$ such that $\left(  u,y\right)  >0$ for all $y\in\Pi$.

Let $\Pi\cap B\left(  Z\right)  =\operatorname{conv}\{u_{1},u_{2}\}$, where
$u_{1}$, $u_{2}\in V$ (recall that $\operatorname{conv}\{A\}$ denotes the
convex hull of $A$). Then
\[
\sum\{y_{i}:y_{i}\in K;\text{ }\left(  y_{i},u\right)  >0\}=\left(
u_{1}+u_{2}\right)  /2.
\]
Hence, $K^{\prime}\subset B(Z)$. Since $B(Z)$ is a closed convex set,
$K^{\prime}=B(Z)$.

$\left(  1\Rightarrow2\right)  $. Let $Z^{\ast}\hookrightarrow B(Z)$ be a
$2$-dimensional subspace. According to the theorem 3, $B(Z)=K^{\prime}$. By
the theorem 4 the incarnating set $K$ for a pair $Z^{\ast}\hookrightarrow
B(Z)$ is situated on beams (or, in other terminology, on semi-axis) of
straight lines that are parallel to edges of $K^{\prime}$. Notice that every
such beam contains just one vector $r$ from $K$. Its length $\left|  r\right|
$ is equal to one half of the length of a edge of $K^{\prime}$ that is
parallel to $r$.

Let $x\in V$ be an exposed point; $\Pi$ be a support line such that $\Pi\cap
B\left(  Z\right)  =\{x\}$. Let $X$ be the one-dimensional subspace of $Z$,
spanned by $x$: $X=\{\lambda x:$ $\lambda\in\mathbb{R\}}$. Let $P:Z\rightarrow
X$ be a projection, parallel to $\Pi$. Let $\Gamma_{x}$ be an arc of the curve
$\partial B\left(  Z\right)  $ from $-x$ up to $x$ in the positive direction.
Certainly, the mapping $P:\Gamma_{x}\rightarrow\left[  -x,x\right]  $ is
one-to-one. By the theorem 4,
\[
\sum\{y_{i}:y_{i}\in K;\text{ }r_{i}\in R_{x}\}=x.
\]

This implies that
\[
-x+\sum\nolimits_{r_{i}\in R_{x}}\left(  x_{2i}-x_{2i-1}\right)  =x.
\]
Hence,
\[
-x+\sum\nolimits_{r_{i}\in R_{x}}P\left(  x_{2i}-x_{2i-1}\right)  =x
\]
The last equality implies that the measure of a set $P\left(  V\cap\Gamma
_{x}\right)  $ on $\left[  -x,x\right]  $ is equal to $0$.

Consider $\Gamma_{x}$ as a graph of a concave function on $\left[
-x,x\right]  $. It is clear that the image of a $0$-measure set under the
mapping $P^{-1}:\left[  -x,x\right]  \rightarrow\Gamma_{x}$ also is of zero
measure. Hence, $V\cup\Gamma_{x}$ is of linear measure $0$.
\end{proof}

\begin{remark}
It is of interest to compare this result with the following of J.
Lindenstrauss $[15]$:

\begin{itemize}
\item  Let $x=\sum a_{i}e_{i}$ and $y=\sum b_{i}e_{i}$ be two elements of
$l_{1}$. Assume that the set
\[
\{a_{i}/b_{i}:i<\infty\}
\]
is dense in the real line. Then the subspace $\operatorname{span}%
\{x,y\}\hookrightarrow l_{1}$ is strictly convex.
\end{itemize}

Using the previous result, a subspace of $l_{1}$ with the same property may be
constructed in other way.

Namely, consider the circle $C=\{(\rho,\varphi):\rho=1\}$ (in polar
coordinates) and put
\[
C_{0}=\{(1,\varphi):\varphi=\pi\sum\nolimits_{k<\infty}\alpha_{k}3^{-k}%
;\alpha_{k}\in\{0,1\}\text{ \ for all }k<\infty\}.
\]
Clearly, $C_{0}$ is the Cantor set on $C$. It will be a set of extreme points
of the unit cell $B(W)$ of the space $W$, if $B(W)$ be defined as
$\operatorname{conv}C_{0}$. By the previous theorem, $W$ is conjugate to a
subspace $W^{\ast}$ of $l_{1}$. Certainly, $\partial W$ is smooth and, hence,
$W^{\ast}$ is strictly convex.
\end{remark}

\section{Properties of spaces $L_{p}\left[  0,1\right]  $ ($1\leq p<\infty$)}

Since the class $\left(  L_{p}\right)  ^{f}$ for $p\in2\mathbb{N}%
\backslash\{2\}$ does not have the amalgamation property, it does not contain
an almost $\omega$-homogeneous space.

\begin{definition}
Let $k<\omega$; $X\in\mathcal{B}$. The space $X$ is said to be almost
$k$-isotropic if for any pair of isometric subspaces $A$, $B$ of $X$ of
dimension $\leq k$ and every $\varepsilon>0$ there exists an isomorphic
automorphism $u:X\rightarrow X$ such that $uA=B$ and $\left\|  u\right\|
\left\|  u^{-1}\right\|  \leq1+\varepsilon$..

$X$ is said to be almost $\omega$-isotropic if $X$ is almost $k$-isotropic for
all $k<\omega$.
\end{definition}

Certainly, every almost $\omega$-isotropic Banach space $X$ is almost $\omega
$-homogeneous with respect to the set $H(X)$ of its finite dimensional
subspaces, and, hence, respect to the set $\frak{M}(X^{f})$.

\begin{remark}
It may be proved that a class $X^{f}$ contains an almost isotropic space
(=almost 1-isotropic) if and only if $\frak{M}(X^{f})$ has the $1$%
-amalgamation property, i.e. if any $V$-formation in $\frak{M}(X^{f})$, which
root is 1-dimensional, is amalgamated in $\frak{M}(X^{f})$. It is known that
all $L_{p}$ ($1\leq p<\infty$) are almost isotropic ($[16]$), i.e. each
$\frak{M}(\left(  L_{p}\right)  ^{f})$ has the 1-amalgamation property.
Notice, that from results of the previous section it follows that
$\frak{M}(\left(  L_{p}\right)  ^{f})$ does not have the 2-amalgamation
property (i.e. for $V$-formations with the 2-dimensional root) for
$p=4,6,8,...$.
\end{remark}

\begin{theorem}
The class $\left(  L_{p}\right)  ^{f}$ for $p\in\lbrack1,\infty]\backslash
\left(  2\mathbb{N}\backslash\{2\}\right)  $ has the amalgamation property.
For such $p$ the space $L_{p}\left[  0,1\right]  $ is almost $\omega$-homogeneous.
\end{theorem}

\begin{proof}
As it was noted by W. Lusky [3], from one W. Rudin's result [4] (in a case of
complex scalars) it follows that the space $L_{p}\left[  0,1\right]  $ is
almost $\omega$-isotropic for $p\neq4,6,8,...$. The same result as in [4] was
obtained by W. Linde [5] in the real case. Hence, $L_{p}\left[  0,1\right]  $
is almost $\omega$-homogeneous with respect to the set $\frak{M}(\left(
L_{p}\right)  ^{f})$. So, the class $\left(  L_{p}\right)  ^{f}$ for
$p\in\lbrack1,\infty]\backslash\left(  2\mathbb{N}\backslash\{2\}\right)  $
has the amalgamation property.
\end{proof}

\begin{corollary}
Let $p\in\lbrack1,\infty)\backslash\left(  2\mathbb{N}\backslash\{2\}\right)
$. For any isometric finite dimensional subspaces $A$, $B$ of $L_{p}\left[
0,1\right]  $ their relative projection constants are equal:
\[
\lambda\left(  A\hookrightarrow L_{p}\left[  0,1\right]  \right)
=\lambda\left(  B\hookrightarrow L_{p}\left[  0,1\right]  \right)  .
\]
\end{corollary}

\begin{proof}
This follows from the preceding theorem and the theorem 8.
\end{proof}

Recall the original Gurarii's definition [2].

\begin{definition}
Let $X$ be a Banach space; $\mathcal{K}$ be a class of Banach spaces. $X$ is
said to be a space of almost universal disposition with respect to
$\mathcal{K}$ if for any pair of spaces $A$, $B$ of $\mathcal{K}$ such that
$A$ is a subspace of $B$ ($A\hookrightarrow B$), every $\varepsilon>0$ and
every isomorphic embedding $i:A\rightarrow X$ there exists an isomorphic
embedding $\hat{\imath}:B\rightarrow X$, which extends $i$ (i.e., $\hat
{\imath}|_{A}=i$) and such, then
\[
\left\|  \hat{\imath}\right\|  \left\|  \hat{\imath}^{-1}\right\|
\leq(1+\varepsilon)\left\|  i\right\|  \left\|  i^{-1}\right\|  .
\]
\end{definition}

Let us show that spaces $L_{p}\left[  0,1\right]  $ ($1\leq p<\infty$) are not
spaces of almost universal disposition with respect to $\frak{M}(\left(
L_{p}\right)  ^{f})$. The proof is different for cases $1<p<2$; $p>2$ (and
$p\neq4,6,8,...$) and for $p=1$.

\begin{definition}
Let $1\leq p\leq\infty$. An operator $u\in B(X,Y)$ is said to be

\begin{itemize}
\item $p$-\textit{absolutely summing,} if there is a constant $\lambda>0$ such
that
\[
(\sum\nolimits_{j<n}\left\|  u\left(  x_{i}\right)  \right\|  ^{p})^{1/p}%
\leq\lambda(\sum\nolimits_{j<n}\left|  \left\langle x_{i},f\right\rangle
\right|  ^{p})^{1/p}%
\]
for any $f\in X^{\ast}$ and any finite set $\{x_{i}:i<n;n<\infty\}\subset X $.
\end{itemize}

Its $p$-\textit{absolutely summing norm }$\pi_{p}\left(  u\right)  $ is the
smallest constant $\lambda$.

\begin{itemize}
\item $p$-\textit{integral,} if there exists a such probability measure $\mu$
and such operators $v\in B(X,L_{\infty}(\mu))$, $w\in B(L_{p}(\mu),Y^{\ast
\ast})$ that $k_{Y}\circ u=w\circ\varphi\circ v$, where $\varphi$ is an
inclusion of $L_{\infty}(\mu)$ into $L_{p}(\mu)$.
\end{itemize}

Its $p$-\textit{integral norm} is given by
\[
\iota_{p}(u)=\inf\{\left\|  v\right\|  \left\|  \varphi\right\|  \left\|
w\right\|  :k_{Y}\circ u=w\circ\varphi\circ v\}.
\]
\end{definition}

\begin{theorem}
Let $X^{f}$ be a class of finite equivalence, which $l_{p}$-spectrum $S\left(
X^{f}\right)  $ is contained in $(1,2]$. Let $W\in X^{f}$ be a space of almost
universal disposition with respect to $\frak{M}\left(  X^{f}\right)  $. Then
for every $r\in S\left(  X^{f}\right)  $ there exists a constant $c_{r}%
<\infty$ such that for every Banach space $Z$ and each finite rank operator
$T:Z\rightarrow W$ its $r$-integral norm $\iota_{r}\left(  T\right)  $ is
estimated by
\[
\iota_{r}\left(  T\right)  \leq c_{r}\pi_{r}\left(  T\right)  .
\]
\end{theorem}

\begin{proof}
Let $T:Z\rightarrow G$ be a finite rank operator with $\pi_{r}\left(
T\right)  =1$. Let $S$ be the unit ball of $Z^{\ast}$, endowed with weak*
(i.e. with $\sigma\left(  Z^{\ast},Z\right)  $-) topology. Let $j:Z\rightarrow
C(S)$ be the canonical embedding: $jz(z^{\prime})=z^{\prime}(z)$ for
$z^{\prime}\in S$. Let $i_{\mu,p}:C(S)\rightarrow L_{p}(S,\mu)$ be the natural
(i.e. identical) embedding; $\mu$ be some measure on $S$. Let $\widetilde
{Z}=jZ\hookrightarrow C(S)$; $i_{\mu,p}(\widetilde{Z})=A_{0}$. The image
$A_{0}$ is a vector subspace of $L_{p}(S,\mu)$. Let $A$ be the closure of
$A_{0}$ in the $L_{p}(S,\mu)$-metric. Let $w_{0}:A_{0}\rightarrow G$ be given
by $w_{0}\circ\left(  i_{\mu,p}\circ j\right)  z=T\left(  z\right)  $.
According to A. Pietsch [17], the measure $\mu$ may be chosen in a such way
that $w_{0}:A_{0}\rightarrow W$ is continuous in $L_{p}(S,\mu)$-metric, and,
hence, admits a unique extension to $w:A\rightarrow W$.

Since $r\in S\left(  X^{f}\right)  $, $L_{p}(S,\mu)$ is finite representable
in $W$ and, because of $T$ is finite rank operator, $A\in\frak{M}\left(
X^{f}\right)  $. Thus, by the definition 11, $w:A\rightarrow W$ may be
extended to a finite rank operator $\widetilde{w}:L_{p}(S,\mu)\rightarrow W$.

So, $T:Z\rightarrow W$ admits a factorization $T=\widetilde{w}\circ i_{\mu
,p}\circ j$. By the definition 12, its $r$-integral norm is uniformly bounded.
\end{proof}

\begin{corollary}
For $1<p<2$ the space $L_{p}$ is not a space of almost universal disposition.
\end{corollary}

\begin{proof}
According to [18] if $1<p<2$ then there exists a Banach space $X$ and a
sequence of finite dimensional operators $u_{n}:X\rightarrow L_{p}$ such that
$\sup\pi_{p}\left(  u_{n}\right)  <\infty$ and $\sup\iota_{p}\left(
u_{n}\right)  =\infty$. Since $p\in S\left(  L_{p}\right)  $, this contradicts
to the previous theorem.
\end{proof}

To consider a case $p>2$ it will be needed a characterization of classes
$X^{f}$ that contain a space of almost universal disposition.

\begin{definition}
A class $X^{f}$ will be called quotient-closed if for every $A\in
\frak{M}\left(  X^{f}\right)  $ and its subspace $B$ the quotient $A/B$
belongs to $\frak{M}\left(  X^{f}\right)  $.
\end{definition}

\begin{theorem}
Every class $X^{f}$ that contains a space of almost universal disposition
(say, $W$) is quotient-closed.
\end{theorem}

\begin{proof}
Let $G<_{f}X$; $F_{0}\hookrightarrow G$; $Z\hookrightarrow G/F_{0}$ be a
finite dimensional space. Let us show that for every $\varepsilon>0$ every
operator $v:Z\rightarrow W$ may be extended to an operator $\widetilde
{v}:G/F_{0}\rightarrow W$ of norm $\left\|  \widetilde{v}\right\|  \leq\left(
1+\varepsilon\right)  \left\|  v\right\|  $. This would imply that for any
pair $A$, $B\in\frak{M}\left(  X^{f}\right)  $, $B\hookrightarrow A$,
$A/B\in\frak{M}\left(  X^{f}\right)  $.

Let $G_{0}\hookrightarrow G$; $V:G\rightarrow E=G/F_{0}$ be a quotient map;
$Z\hookrightarrow E$ and $v:Z\rightarrow W$. Let $T:G\rightarrow W$. Put
$F=V^{-1}\left(  Z\right)  $. Then the operator $T\circ v\circ\left(
V\mid_{F}\right)  $ may be extended to $w:G\rightarrow W$ with a norm
$\left\|  w\right\|  \leq\left(  1+\varepsilon\right)  \left\|  v\right\|  $.

Since $V^{-1}\left(  0\right)  \subset F$, $w\left(  V^{-1}\left(  0\right)
\right)  =0$ and, hence, $w=v_{1}\circ V$, where $v_{1}:E\rightarrow W$;
$\left\|  v_{1}\right\|  \leq\left(  1+\varepsilon\right)  \left\|  v\right\|
$. Certainly, $v_{1}$ is the desired extension of $v$.
\end{proof}

\begin{corollary}
For $2<p<\infty$ the space $L_{p}$ is not a space of almost universal disposition.
\end{corollary}

\begin{proof}
It is obvious that for every quotient-closed class $X^{f}$ its $l_{p}%
$-spectrum is:
\[
S\left(  X^{f}\right)  =\left[  \inf S\left(  X^{f}\right)  ;\sup S\left(
X^{f}\right)  \right]  .
\]
Since $S(L_{p})=\{2,p\}$ for $p\in\left(  2,\infty\right)  $, the class
$\left(  L_{p}\right)  ^{f}$ cannot be quotient closed.
\end{proof}

It is remain the case $p=1$.

\begin{theorem}
$L_{1}\left[  0,1\right]  $ is not the space of almost universal disposition.
\end{theorem}

\begin{proof}
If the class $X^{f}$ is quotient closed and $l_{1}$ is finitely representable
in $X$ then $X^{f}=\left(  l_{\infty}\right)  ^{f}$. So, the previous theorem
yields the needed result. Another proof follows from [6], prop. 1.11.
\end{proof}

\section{References}

\begin{enumerate}
\item  J\'{o}nsson B. \textit{Universal relational systems}, Math. Scand.
\textbf{4} (1956) 193-208

\item  Gurarii V.I.\textit{\ Spaces of universal disposition, isotropic spaces
and the Mazur problem on rotations in Banach spaces}, Sibirsk. Mat. Journ. (in
Russian) \textbf{7} (1966) 1002-1013

\item  Lusky W. \textit{Some consequences of W. Rudin's paper ''}$L_{p}%
$\textit{\ - isometries and equimeasurability''}, Indiana Univ. Math. J.
\textbf{27:5} (1978) 859-866

\item  Rudin W. $L_{p}$\textit{-isometries and equimeasurability,} Indiana
Math. J. \textbf{25:3} (1976) 215-228

\item  Linde W. \textit{Moments of measures in Banach spaces}, Math. Ann.
\textbf{258} (1982) 277-287

\item  Pisier G. \textit{Counterexample to a conjecture of Grothendieck}, Acta
Math. \textbf{151} (1983) 181-208

\item  Positselskii E.D., Tokarev E.V. \textit{Amalgamation of classes of
Banach spaces }(in Russian), \textbf{XI} All-Union School on Operator Theory
in Functional Spaces. A book of abstracts, Tcherlyabinsk (1986) 89

\item  Zalgaller V.A., Reshetniak Yu.G. \textit{On straightened curves,
additive vector-functions and merging of segments }(in Russian), Vestnik LGU
\textbf{2} (1954) 45-67

\item  Stern J. \textit{Ultrapowers and local properties in Banach spaces},
Trans. AMS \textbf{240} (1978) 231-252

\item  Tokarev E.V. \textit{Injective Banach spaces in the finite equivalence
classes }(transl. from Russian), Ukrainian Mathematical Journal \textbf{39:6}
(1987) 614-619

\item  Lindenstrauss J., Pe\l czy\'{n}ski A. \textit{Absolutely summing
operators in }$\mathcal{L}_{p}$\textit{-spaces and their applications}, Studia
Math. \textbf{29} (1968) 275-326

\item  Garling D.J.H., Gordon Y. \textit{Relations between some constants
associated with finite dimensional Banach spaces}, Israel J. Math.
\textbf{9:3} (1971) 346-361

\item  Kadec M.I., Mitiagin B.S. \textit{Complemented subspaces in Banach
spaces}, Uspechi Math. Nauk \textbf{28:6 }(1973) 77-94

\item  Straszewicz S. \textit{\"{U}ber exponierte Punkte abgeschlossener
Punktmengen}, Fund. Math. \textbf{24} (1935) 139-143

\item  Lindenstrauss J. \textit{On the extension of operators with a finite
dimensional range}, Illinois J. Math. 8 (1964) 488-499

\item  Rolewicz S. \textit{Metric linear spaces}, Warszawa: Monografie
Matematyczne \textbf{56} (1972)

\item  Pietsch A. \textit{Absolut }$p$\textit{-summierende Abbildungen in
normierten R\"{a}umen}, Studia Math. \textbf{28} (1967) 333-353

\item  Pe\l czy\'{n}ski A. $p$\textit{-integral operators commuting with group
representations and examples of quasi }$p$\textit{-integral operators, which
are not }$p$\textit{-integral}, Studia Math. \textbf{33} (1969) 63-70
\end{enumerate}
\end{document}